
\input amssym
\magnification = \magstep1
\centerline{{\bf Mertens Sums requiring Fewer Values of the M\"obius
Function}}
\medskip 
\centerline{{\it (dedicated to the memory of Yu. V. Linnik)}}
\bigskip
\centerline{{\it M. N. Huxley and N. Watt}}
\bigskip
\centerline{{\bf 1. Introduction}}
\medskip
The sieve of Eratosthenes will find the prime numbers in
$N+1,\dots ,N^2$ provided that we know all the primes in $2,\dots ,N$.
In particular the sieve gives a relation for the function $\pi (x)$ that
counts the number of primes less than or equal to~$x$:
$$
\pi (N^2)  =   \pi (N) -1 +
\mathop{\sum _{d \leq N^2}}_{P(d)\leq N}\mu (d)
\left[N^2\over d\right],
\eqno (1.1)
$$
where $\mu (d)$ is the M\"obius function (which is $(-1)^\nu $ when
$d$ has $\nu $~prime factors, all different, but 0 when $d$ has any
prime factor repeated),
while $P(d)$ is the greatest non-composite divisor of $d$, and
$[x]=\max\{ m\in{\Bbb Z} \,:\, m\leq x\}$. The numbers~$d$ in~(1.1) are
constructed as
products of the known primes in $2,\dots ,N$, so the values~$\mu (d)$
can be read off. In general, given a number~$n$, it is very difficult
to factorise~$n$ and so find~$\mu (n)$. Thus the Mertens sum
$$
M(x)  =  \sum _{n \leq x}\mu (n)
\eqno (1.2)
$$
is difficult to calculate from the definition. The Dirichlet series
$\sum \mu (n)/n^s$ is~$1/\zeta (s)\,$ (the reciprocal of the Riemann
zeta function), and, according to folklore, the fastest method of
calculating~$M(x)$ is by Perron's contour integral formula for the sum
of the coefficients of a Dirichlet series.

In this paper we discuss a family of identities which allow $M(N^d)$ to
be calculated for each positive integer $d$ as a sum of no more than
$O_d \left( N^d(\log N)^{2d - 2}\right)$ terms, each a product of the
form $\mu(n_1) \cdots \mu(n_r)$ with $r\leq d$ and
$\{ n_1, \ldots , n_r\}\subseteq\{1,\ldots ,N\}$. In Theorem~1, below,
we state a more complicated form of these identities, in which
each of the variables of summation $n_j\,$ ($j=1,\ldots ,r$) can have
its own  independent range of summation:  $1,\ldots , N_j\,$ (say).

We actually treat the more general M\"obius sum
$$
M(g , x) = \sum _{n \leq x}\mu(n)g(n),
\eqno(1.3)
$$
where $g(n)$ can be any totally multiplicative arithmetic function,
that is, $g(rs) = g(r)g(s)$ holds for any positive integers $r$
and~$s$. The relevant identity when $d = 1$ is (of course) the
definition~(1.3). The case $d = 2$ is the next simplest. Let
${\bf m}(g,N)$ be the column-matrix
$(\mu (1)g(1),\ldots ,\mu (N)g(N))^{\rm T}$,
and let $A(g, N)$ be the $N\times N$~matrix with elements
$$
a_{mn}(g , N) = \sum _{k \leq {\scriptstyle N^2\over\scriptstyle mn}}g(k)
\qquad\hbox{($m,n\in\{ 1,\ldots ,N\}$).}
\eqno(1.4)
$$
Then
$$
M(g , N^2) = 2M(g , N) - \left( {\bf m}(g,N)\right)^{\rm T} A(g , N)
    {\bf m}(g , N)\;.
\eqno (1.5)
$$
In the general case, when $d,K,N\in{\Bbb N}$ satisfy $d\geq 2$ and
$K\geq N> K^{1/d} -1$, we have:
$$
\eqalign{
M(g , K)  &=  dM(g , N) \cr
&\quad - \sum _{r = 2}^d (-1)^r \,\null_dC_r
\mathop{\sum _{n_1 \leq N}\dots \sum _{n_r \leq N}
\sum _{k_1}\dots \sum _{k_{r-1}}}
_{n_1n_2\dots n_rk_1k_2\dots k_{r-1} \leq K}
g(k_1 \cdots k_{r-1}) \prod_{i=1}^r \mu (n_i)g(n_i)\;,
\cr}
\eqno(1.6)
$$
where $\null_d C_r = d(d-1)\cdots (d-(r-1))/(r!)$.

Note that (1.5) is just the special case $d=2$, $K=N^2$ of (1.6).
Moreover, (1.6) is itself a special case of another identity
(that stated in Theorem~1, below), in which
the single range of summation $1,\ldots ,N$ is replaced by
$d$ independent ranges of summation.
In order to state this more general identity we require some more notation.
\par
Let $d$ be a positive integer greater than $1$. Let $V = v_1v_2\dots
v_d$ be a word of length $d$ in the alphabet $\{0,1\}$.
The support of a word~$V$ is the set of indices~$i$ for which $v_i = 1$.
The weight~$w(V)$ of a word~$V$ is the size of the support, so that
$w(V) = \sum v_i$. The combinatorial M\"obius function, which we
write as~$\mu ^*$ to distinguish it from the number-theoretic
function~$\mu $, is $\mu ^*(V) = (-1)^{w(V)}$.

Let $N_1,\ldots , N_d$ be positive integers.
For each word~$V$, and each $L\in{\Bbb N}$,
let the notatation $\sum _1^L (V)$ signify summation over
$n_1,\dots ,n_d$ in the ranges $n_i = 1,\dots ,L$ when $v_i = 0$, but
$n_i = 1,\dots ,N_i$ when $v_i = 1$.
When $L = 1$ and $v_i=0$, the variable of summation $n_i$ effectively
becomes `frozen', meaning that its range of summation is then just the
single-element set~$\{ 1\}$.

Let $K$ be a positive integer that is less than $(1+N_1)(1+N_2)\dots (1+N_d)$.
If $n_1,\ldots , n_d$ are integers satisfying the condition
$n_1n_2\dots n_d \leq K$, then $n_i \leq N_i$
holds for at least one index~$i$. It therefore follows by the
inclusion-exclusion principle of combinatorics that if
$f: {\Bbb N}^d\rightarrow{\Bbb C}$ is such that one has
$|f(n_1,\ldots , n_d)|>0$ only when $n_1 n_2 \cdots n_d \leq K$, then
$$
\sum _1^K(00\dots 0) f( n_1 , \ldots , n_d)
=\sum_{r=1}^d (-1)^{r-1} \sum_{V\,:\,\omega(V)=r}
\ \sum_1^K (V) f( n_1 , \ldots , n_d)\;,
\eqno (1.7)
$$
or, to put it more elegantly,
$\sum _V\mu ^*(V)\sum _1^K(V) f( n_1 , \ldots , n_d) = 0$.

\proclaim Theorem 1. When $g(n)$ is a totally multiplicative arithmetic
function, and $d$, $N_1,\ldots , N_d$ and $K$ are as above, we have:
$$
\eqalign{
M(g,K) &= \sum_{i=1}^d M\!\left( g,\min\{ N_i , K\} \right) \cr
&\quad - \sum _{V\, :\, w(V) \geq 2}(-1)^{w(V)}
   \sum _1^1(V) \!\!\mathop{\sum _{k_1}\dots
   \sum _{k_{\omega(V)-1}}}_{k_1\dots k_{\omega(V)-1}
\leq {\scriptstyle K\over\scriptstyle n_1\dots n_d}} \!\!\!\!\!\!
g(k_1 \cdots k_{\omega(V) -1}) \prod_{i=1}^{d} \mu (n_i)g(n_i)\;.
\cr}
\eqno(1.8)
$$

\noindent
{\bf Proof.}\quad
We apply (1.7) with $f$ given by:
$$
f(n_1 , \ldots , n_d) =\mathop{\sum _{k_1}\dots \sum _{k_{d-1}}}
   _{k_1\dots k_{d-1} \leq {\scriptstyle K\over\scriptstyle n_1\dots n_d}}
   \mu(n_1)\dots \mu(n_d)g(k_1\dots k_{d-1}n_1\dots n_d).
\eqno (1.9)
$$
For the word $V =11\dots 1$ with $w(V) = d$, we have
$$
\sum _1^K(11\dots 1) = \sum _1^1(11\dots 1).
$$
All other words~$V$ have $v_j = 0$ for at least one index~$j$, so the
corresponding summand $n_j$ runs over the full range from~1 to~$K$.
For these words $V$ we carry out the following `contraction step'.
Take an index~$j$ for which $v_j = 0$. We sum over
$n_j$ and~$k_{d-1}$ first, observing that by M\"obius inversion we have:
$$
\eqalign{
\sum _{n_j = 1}^K \,
&\!\!\!\!\!\!\!\!\!\!\!\!\!\!\!\!\!\!\!\!\!\!\!\!
\sum _{\qquad\qquad\ k_{d-1}\leq
  {\scriptstyle K\over\scriptstyle n_1\dots n_dk_1\dots k_{d-2}}}
   \mu (n_j)g(n_jk_{d-1}) \cr
&= \!\!\!\!\!\!\!\!\!\!\!\!\!\!\!\!\!\!\!\!\!\!\!\!\!\!\!\!\!\!\!\!\!
\!\!\!\!\!\!\!\!\!\!\!\!\!\!\!
\!\!\!\!\!\!\!\sum _{\qquad\qquad\qquad\qquad\quad\
m \leq {\scriptstyle K\over\scriptstyle n_1\dots n_{j-1}n_{j+1}
   \dots n_dk_1\dots k_{d-2}}}g(m)  \sum _{n_j|m}\mu (n_j)\cr
    &= \cases{
    g(1) = \mu (1)g(1)
    &if $n_1\dots n_{j-1}n_{j+1}\dots n_dk_1\dots k_{d-2}\leq K$, \cr
    0 &otherwise.}
\cr}
$$
We thereby find that the value of the relevant sum
over $n_1, \ldots , n_d$ and $k_1 , \ldots , k_{d-1}$
is unchanged when we
omit~$k_{d-1}$ and freeze~$n_j$ as the fixed value $n_j = 1$.
\par
We repeat the contraction step for every index~$j$ with $v_j = 0$,
freezing the corresponding variable as $n_j = 1$, and removing the last
variable~$k_i$. Exceptionally, when $V$~is $00\dots 0$, we can remove
$k_{d-1},k_{d-2},\dots ,k_1$, and freeze $n_d,n_{d-1},\dots ,n_2$, but
the sum over~$n_1$ remains over the range $1,\dots ,K$, giving the
term~$M(g,K)$ on the left of~(1.8). The summation identity (1.7), when applied
with $f$ given by (1.9), contracts to give (1.8).\hfill$\blacksquare$

\medskip

In (1.5), (1.6) and Theorem~1, we require the total multiplicativity of~$g$
only in order to `separate variables' (as, in (1.6) for example, we separate
$k_1,\ldots ,k_{r-1}$ from $n_1,\ldots ,n_r$ by means of the identity
$g(k_1 \cdots k_{r-1}n_1 \cdots ,n_r)
= g(k_1 \cdots k_{r-1})g(n_1) \cdots g(n_r)$). Indeed, (1.8) gives a
formula for the M\"{o}bius function itself, for we can apply (1.8) to
each term in the difference $M(g,K) - M(g,K-1) = \mu(K)g(K)$, and we
can then divide through by~$g(K)$ to obtain a formula for $\mu(K)$ that
is independent of~$g$.
This formula for $\mu(K)$ may also be deduced from the identity
$$
{1\over \zeta(s)}\prod_{j=1}^d \biggl( 1- \zeta(s)\sum_{n=1}^{N_j}
   {\mu(n)\over n^s}\biggr)
=\zeta^{d-1}(s)\prod_{j=1}^d
\,\sum_{n = 1+N_j}^{\infty} {\mu(n)\over n^s}\qquad\hbox{(${\rm Re}(s) > 1$),}
\eqno(1.10)
$$
through multiplying out the brackets on the left-hand side, and then
computing the coefficient of $n^{-Ks}$ on each side of the resulting
identity, subject to the hypothesis that the product
$(1+N_1)\cdots (1+N_d)$ be greater than $K$.
This approach yields a second proof of Theorem~1. We prefer the first
proof due to its more obvious connection with Meissel's
identity~[{\bf 8} p~303],

$$\sum_{n\leq x} \left[ {x\over n}\right] \mu(n)
=\cases{1 &if $x\geq 1$, \cr 0 &if $1>x>0$, \cr}\eqno(1.11)$$
which was the initial source of inspiration for our work.
\par
Given any $K\in{\Bbb N}$, any integer $d\geq 2$,
and any $\theta_1,\ldots ,\theta_d >0$ with $\theta_1 + \cdots +\theta_d =1$,
it follows from Theorem~1 that (1.8) will hold when one has also
$N_j=[K^{\theta_j}]$,  for $j=1,\ldots ,d$. Theorem~1 therefore offers
considerably more flexibility of application than (1.6) does. Although
we believe Theorem~1 to be new (in respect of the flexibility in the
choice of $N_1,\ldots , N_d$), the special cases of it that are displayed in
(1.5) and (1.6) are known results. The result (1.5) is
contained in Vaughan's (slightly more complicated) identity
[{\bf 13} equation~(18)]
(essentially the special case when $u=\sqrt{X}$, and so $S_3=0$), and one can
find in equation (13.38) of~[{\bf 5}],
for example, a formula for~$\mu(n)$ that is equivalent to what we have
in~(1.6). It is, moreover, clear that even our identity in~(1.8) is
akin to formulae of Heath-Brown for sums involving  $\Lambda(n)$, the
von~Mangoldt function:
compare (1.10), from which (1.8) may be deduced,
with Lemma~1 of~[{\bf 2}].
The earliest formula of this type is due to
Linnik himself in~[{\bf 6},{\bf 7}].
\par
We shall refer to the case of~(1.3) (or of (1.4), (1.5), (1.6), or~(1.8)),
where the function~$g(n)$ takes the constant value~$1$, as the
principal case. The main focus of our work has been on the principal
case of the identity~(1.5). Indeed, all subsequent sections of this
paper are exclusively devoted to matters connected with this single topic,
such as (for example)
questions concerning certain properties of the $N\times N$ matrix
$A=A(N)$ that occurs in the principal case of (1.5) and has, by (1.4),
elements $a_{mn}=[N^2 /(mn)]\in{\Bbb N}$.
In Section~2 we discuss matters related to the spectral
decomposition of $A=A(N)$.
In the third (and final) section we discuss decompositions (spectral
and otherwise) of the quadratic form ${\bf m}^T A {\bf m}$,
where ${\bf m}={\bf m}(N)$ is the column-matrix
$(\mu(1),\ldots ,\mu(N))^{\rm T}$
that occurs in the principal case of (1.5).
\par
We consider especially the principal case of (1.5), in the hope that it
(modified as necessary) might lead to a new proof of the prime
number theorem, or even some new upper bound for the Mertens sum~$|M(x)|$.
The following parts of this paper report what we have discovered in
the search for such an application of~(1.5).
\par
One of our findings is that the
matrix~$A(N)$, which (clearly) is real and symmetric,
has one exceptionally large positive eigenvalue,
approximately~$N^2\zeta (2)$,
with eigenvector approximately~$(1,1/2,1/3,\ldots ,1/N)^{\rm T}$.
Calculations by the second author show that the
second-largest eigenvalue of $A(N)$ lies in an interval of the form
$[d_4 N + o(N) , c_4 N + o(N)]$, where $c_4$ and $d_4$ are constants
that are approximately $-0.496$ and $-0.572$, respectively: for more details,
see (2.7), (2.14), (2.20) and (2.21) below. Hence, for $N$ sufficiently large,
the quadratic form on the right-hand side of (1.5) is neither positive
definite nor negative definite in the principal case.
\par
By the principal case of (1.6), we have a sequence of formulae
through which each of $M(N^2),M(N^3),M(N^4),\ldots$ is expressed in
terms of $\mu(1),\ldots ,\mu(N)$. Although the first of these formulae,
the principal case of~(1.5), may be considered analogous to the sieve of
Eratosthenes~(1.1), there seems to be no version of~(1.1)
for~$\pi (N^3)$,  because unwanted numbers of the form~$pq$, where $p$
and~$q$ are both primes greater than~$N$, survive the sieve process
(``Gnoggensplatts'' in Greaves's lectures on {\it Sieve Methods}).
\par
A connection between Mertens sums and certain symmetric matrices $U_n\,$
($n\in{\Bbb N}$),
that bear some resemblance to our matrices $A(N)\,$ ($N\in{\Bbb N}$)
has previously been established by
Cardinal~[{\bf 1}].
To define Cardinal's matrix $U_n$, one first takes
$\sigma_1 <\sigma_2 <\ \cdots\ <\sigma_s$ to be the elements of the set
${\cal S}={\cal R}\cup\{ [n/\rho]\,:\,\rho\in{\cal R}\}$, where
${\cal R}=\{ \rho\in{\Bbb N}\,:\, \rho \leq\sqrt{n}\}\,$
(it follows that $0\leq 2[\sqrt{n}\,] - s\leq 1$). Then $U_n$ is the
$s\times s$ matrix with elements $u_{ij}=[n/(\sigma_i \sigma_j)]$.
In Propositions 21 and~22 of~[{\bf 1}],
it is shown that one has $T_n U_n^{-1} T_n = V_n$, where $T_n$ and
$V_n$ are the $s\times s$ matrices with elements
$t_{ij}=|[2,s+1]\cap\{i+j\}|$ and $v_{ij}=M(u_{ij})$, respectively.
\par
In the cases where $n$ is a perfect square, so that $n = N^2$ for some
integer~$N$, then $|{\cal R}| = N$, and the $N\times N$ principal
submatrix of~$U_n$ consisting of the array of elements from the first
$N$ rows and first $N$ columns of~$U_n$ is our matrix $A(N)$: since
$2N-1\leq s\leq 2N$, we can say that $A(N)$ constitutes (exactly, or
approximately) the top left-hand quarter of Cardinal's matrix~$U_n$.
In these same cases, Cardinal's identity
$T_n U_n^{-1} T_n = V_n$ implies that $v_{11}$, which is $M(N^2)$,
will be equal to the sum of all $s^2$ of the elements of the inverse of
the matrix $U_n = U_{N^2}$: we obtain a formula for $M(N^2)$ thereby
that seems quite different from what we see in the principal case
of~(1.5).
\par
As Cardinal
observes in Theorem~24 and Remark~25 of~[{\bf 1}],
information about small eigenvalues of the matrix
$V_n^{-1}=T_n^{-1} U_n T_n^{-1}$ might lead to new upper bounds
on~$M(x)$. In this respect, the connection that we have found between
$M(x)$ and~$A(N)$ is quite different from Cardinal's connection between
$M(x)$ and~$U_n$, for it is the larger eigenvalues of~$A(N)$ and their
eigenvectors that matter most in the principal case of~(1.5): see, for
example, equation~(3.3), below.
\par
We have scarcely considered non-principal cases of (1.5), (1.6),
or~(1.8). Certain non-principal cases of~(1.5) may merit further
investigation.
The first case is when $g(n) = \chi (n)$, a non-principal Dirichlet
character to some modulus~$q > 1$. The sums
$\sum _{\ell \leq x}\chi (\ell)$ that we use to construct the matrix
elements~$a_{mn}(\chi , N)$ in~(1.4) are periodic step functions
of~$x$, whose period is~$q$ or some proper factor of~$q$.
In contrast to the principal case, when the set of elements of the
matrix $A(N)$ in~(1.5)) contains at least~$N$ dirrerent integers,
namely $[N^2 /1],[N^2 /2],\dots ,[N^2 /N]$, there is a single finite set,
$\{ \sum_{0<\ell\leq L}\chi(\ell) :  L\in (0,q]\cap{\Bbb Z}\}$, that
contains all the elements of all the matrices
$A(\chi , 1), A(\chi , 2) , A(\chi ,3) , \ldots\ $.
For $\chi$ real, $A(\chi , N)$ will, of course, be real and symmetric
just like~$A(N)$.
\par
A case of~(1.3) known to be related to the prime number theorem is
when $g(n) = 1/n\,$
(see page~248 of~[{\bf 9}],
for example).
More generally, when $g(n) = n^{-s}$ for some fixed complex number~$s$,
then the sum~$M(g, x)$ in~(1.3) becomes a partial sum for the Dirichlet
series for~$1/\zeta (s)$.
If, for some $\sigma_{0}\in [1/2 , 1)$,
the only zeros of $\zeta (s)$ with real parts greater than $\sigma_{0}$
are a pair of simple zeros, $\rho$ and~$\overline{\rho}\,$ (say),
and if we put $g(n) = n^{-\rho}\,$ ($n\in{\Bbb N}$), then the
sum~$M(g,x)$ in~(1.3) will grow logarithmically in~$x$.
\par
Another interesting case of (1.3) to (1.5) is when $g(n) = \lambda (n)$,
the Liouville function, which
is the projection of the M\"obius function~$\mu $ onto the space of
totally multiplicative arithmetic functions. In this case $M(g,x)$ grows
like~$x/\zeta (2)$.

\beginsection 2. Elementary Estimates for Eigenvalues and an Eigenvector

\medskip

Let $N$ be a given positive integer. Since the matrix $A = A(N)$, in
the principal case of~(1.5), is both real and symmetric, it has
eigenvalues $\lambda_1 \leq \lambda_2 \leq \  \ldots\ \leq\lambda_N$
with corresponding eigen-(column-)vectors of unit length
${\bf e}_1 , \ldots , {\bf e}_N$ that form an orthonormal basis
of~${\Bbb R}^N$. When ${\bf v}\in{\Bbb R}^N$, one has
$$
{\bf v}^T A {\bf v} =
\sum_{k=1}^N \lambda_k \left( {\bf e}_k \cdot {\bf v}\right)^2
\eqno(2.1)
$$
as a consequence of the spectral decomposition
$A = \sum_{k=1}^N \lambda_k {\bf e}_k {\bf e}_k^T$, and Parseval's
identity gives
$$
\sum_{k=1}^N \left( {\bf e}_k \cdot {\bf v}\right)^2 = {\bf v} \cdot {\bf v}
=  \| {\bf v} \|^2 .
\eqno(2.2)
$$
In order to study the terms appearing in (2.1) and~(2.2), we estimate:
\smallskip
(a)\qquad ${\rm Tr}(A) = \sum a_{nn}\,$ (the trace of the matrix~$A$),
\smallskip
(b)\qquad ${\rm Tr}(A^2) = {\rm Tr}(A^T A)=\sum \sum a_{mn}^2$,
\smallskip
(c)\qquad ${\bf f}^T A {\bf f}$,
where ${\bf f} = (1,{1\over 2},{1\over 3},\ldots ,{1\over N})^T$,
\smallskip
(d)\qquad ${\bf w}^T A {\bf w}$,
where ${\bf w} = {\bf u} - \| {\bf f}\|^{-2} ({\bf f}\cdot {\bf u}) {\bf f}$,
with ${\bf u} = (1,1,\ldots ,1)^T \in{\Bbb R}^N$.

\medskip

We use the following notation:
$$
\zeta_j=\sum_{m=1}^N m^{-j}\;, \quad \delta
=\sum_{m\leq N}\sum_{n\leq N}{\left\{ N^2 /(mn)\right\}\over mn}\quad
{\rm and}\quad
\phi =
{1\over N^2}\sum_{m\leq N}\sum_{n\leq N} \left\{ {N^2 \over mn}\right\}^2\;,
$$
where $\{ t\} = t - [t]$ (the fractional part of~$t$).
Taking (b) first, we simply observe that
$$
{\rm Tr}\left( A^2\right)
= \sum_{m\leq N} \sum_{n\leq N} \left[ {N^2 \over mn}\right]^2
= \sum_{m\leq N} \sum_{n\leq N} \left( {N^2 \over mn}
- \left\{ {N^2 \over mn}\right\}\right)^2
=\zeta_2^2 N^4 + \left(\phi - 2 \delta\right) N^2\;.
\eqno(2.3)
$$
Since ${\rm Tr}(A^2) = \lambda_1^2 +\ \cdots\ + \lambda_N^2$, and since
$\delta\geq 0$ and $\phi < 1$,
the identity (2.3) shows already that
$\lambda_N < \zeta_2 N^2 + (2\zeta_2)^{-1}$.

\par
Regarding~(c), we are content to note that
$$
{\bf f}^T A {\bf f} =
\sum_{m\leq N}\sum_{n\leq N} {\left[ N^2 /(mn)\right]\over mn}
=   \sum_{m\leq N}\sum_{n\leq N}
{\left( {N^2 /(mn)} -\left\{ {N^2 /(mn)}\right\}\right)\over mn}
=   \zeta_2^2 N^2 - \delta\;.
\eqno(2.4)
$$
We have here $\| {\bf f}\|^2 = \zeta_2$, so by Rayleigh's Principle it
follows from~(2.4) that
$$
\zeta_2 N^2 - {\delta\over\zeta_2} \leq \lambda_N\;.
\eqno(2.5)
$$
By (2.5) and the point noted immediately below~(2.3), we conclude that
$$
- {(1 + \log N)^2\over\zeta_2} < \lambda_N - \zeta_2 N^2
< {1\over 2\zeta_2}\;.
\eqno(2.6)
$$
As
$0\leq \delta < \zeta_1^2\leq \zeta_0\zeta_2 = N\zeta_2\leq N^2 \zeta_2^2$,
the lower bound on~$\lambda_N$ in~(2.5) is non-negative, and so we may
deduce from it that
$\lambda_N^2 \geq (\zeta_2 N^2 - \delta \zeta_2^{-1})^2
= \zeta_2^2 N^4 - 2\delta N^2 + \delta^2 \zeta_2^{-2}$:
this, together with the evaluation of~${\rm Tr}(A^2)$ in~(2.3), is enough
to show that
$$
\lambda_1^2 + \ \cdots\ + \lambda_{N-1}^2
\leq \phi N^2 - \delta^2 \zeta_2^{-2} < N^2\;.
\eqno(2.7)
$$
From the way we have ordered the eigenvalues, the bound (2.7) implies:
$$
\left| \lambda_k\right| < {N\over\sqrt{\min\{ k , N-k\}}}\qquad
\hbox{($k=1,2,\ldots , N-1$).}
\eqno(2.8)
$$

\par
In view of (2.6) and~(2.8), it is clear that for $N$ large,
$\lambda_N$ will be exceptionally large, compared with all other
eigenvalues of~$A$.
Accordingly we consider first the corresponding eigenvector~${\bf e}_N$,
before discussing the estimation (a) of~${\rm Tr}(A)$.
Putting $F_N = {\bf e}_N \cdot \hat{\bf f}$, where
$\hat{\bf f} = \|{\bf f}\|^{-1} {\bf f}$,
we find by (2.4) and~(2.6), and (2.1), (2.8) and~(2.2), that
$$
\lambda_N - \left( \textstyle{1\over 2} + (1 + \log N)^2\right)
< \hat{\bf f}^T A \hat{\bf f}
< \lambda_N F_N^2 + N\left( 1 - F_N^2\right)\;.
$$
For $N > 1$ we have $\lambda_N > N$ (this follows by~(2.6) when $N\geq 3$),
and so, by comparison of the upper and lower bounds for
$\hat{\bf f}^T A \hat{\bf f}$ that were just obtained, we deduce that
$$
1  \geq F_N^2 > 1 - {\left( \textstyle{1\over 2} + (1 + \log N)^2\right)
\over \left( \lambda_N - N\right)}\;.
$$
Choosing the $\pm$-sign so that ${\pm}F_N = |F_N|$, we therefore find
from~(2.6) that
$$
\left\| {\bf e}_N - \left( \pm\hat{\bf f} \right)\right\|
= \sqrt{2\left( 1 - |F_N|\right)}
= \sqrt{ {2\left( 1 - F_N^2\right)\over 1 + |F_N|} }
= O\left( {\log N\over N}\right) \;.
\eqno(2.9)
$$

\par
We now come to the task mentioned in~(a) above, which is the estimation
of the sum $S = {\rm Tr}(A) = \sum a_{n n}$. We pick a positive
integer~$K$, and we divide the original sum~$S$ into two parts: $S_1$,
which has the terms with $n^2 \leq N^2 /(K + 1)$,
and $S_2$, which has the terms with $N^2\geq n^2 > N^2 /(K + 1)$
(so that $a_{n n} =[N^2 /n^2] = k$ for some $k\in\{ 1,\dots ,K\}$). We have
$$
\eqalign{
S_1 = \sum _{n^2 \leq N^2/(K + 1)}a_{nn}
 &= \sum _{n \leq  N/\sqrt{K + 1}}\left({N^2\over n^2}+O(1)\right) \cr
 &= N^2\left(\zeta_2 - \int _{N/\sqrt{K + 1}}^N x^{-2}\,{\rm d}x
   + O\left({K\over N^2}\right)\right) + O\left({N\over \sqrt{K}}\right) \cr
 &=\zeta_2 N^2-N\sqrt{K}+N+O\left(K + {N\over \sqrt{K}}\right).
\cr}
$$
The sum~$S_2$ is more complicated. We have
$$
\eqalign{
S_2 =\sum _{k= 1}^K \sum _{{N\over\sqrt{k+1}} < n \leq {N\over\sqrt{k}}} k
 &=\ \sum \sum _{\!\!\!\!\!\!\!\!\!\!\!\!\!\!{1\leq\ell\leq k\leq K}}
\left( \left[ {N\over \sqrt{k}}\right]
 - \left[ {N\over \sqrt{k+1}}\right]\right)  \cr
 &= \sum _{\ell = 1}^K \left( \left[ {N\over\sqrt{\ell}}\right]
 - \left[ {N\over\sqrt{K+1}}\right]\right)
= \sum _{\ell = 1}^K {N\over \sqrt{\ell}} -{K N\over \sqrt{K+1}}+O(K).
\cr}
$$
Let
$$
g(\ell ) = 2\sqrt{\ell }-2\sqrt{\ell -1}-{1\over \sqrt{\ell }}
   =  {1\over \sqrt{\ell }(\sqrt{\ell }+\sqrt{\ell -1})^2}
   \quad\hbox{($\ell\in{\Bbb N}$)} \quad\ {\rm and}\quad\
\alpha = \sum _{\ell = 1}^\infty g(\ell ).
$$
Then
$$
\sum _{\ell = 1}^K {1\over \sqrt{\ell }}
   = \sum _{\ell = 1}^K\left(2\sqrt{\ell }-2\sqrt{\ell -1}-g(\ell )\right)
   = 2\sqrt{K}-\alpha +O\left({1\over \sqrt{K}}\right).
$$
Hence
$$
S_2 = 2N\sqrt{K}-\alpha N-{NK\over \sqrt{K+1}}
   +O\left({N\over \sqrt{K}} + K\right)
   = N\sqrt{K}-\alpha N+O\left({N\over \sqrt{K}} + K\right),
$$
and so, putting $K=[N^{2/3}]$, we get:
$$
{\rm Tr}(A) = S_1 + S_2
= \zeta_2 N^2 - (\alpha -1)N + O\left(N^{2/3}\right)\;.
\eqno(2.10)
$$
\par
By (2.10) and~(2.6), it follows that
$$
\lambda_1 + \ \cdots\ + \lambda_{N-1}
= -(\alpha - 1) N + O\left( N^{2/3}\right) \;.
\eqno(2.11)
$$
By equations (1.11) to (1.13) of~[{\bf 4}]
and the case $K=1$ of
of equation (B.24) of~[{\bf 9}]
(itself an application of the Euler-Maclaurin summation formula), we find
that for $\sigma\in (0,1)\cup(1,\infty)$ and $K\in{\Bbb N}$,
$$
\eqalignno{
\sum_{\ell = 1}^K {1\over \ell^{\sigma}}
  &= {K^{1-\sigma} \over 1 - \sigma} + \zeta(\sigma)
+ {\theta(K , \sigma)\over K^{\sigma}} &(2.12)\cr
 &= {\theta(K , \sigma)\over K^{\sigma}} + {K^{1-\sigma} - 1 \over 1 - \sigma}
+ \gamma + \sum_{j=1}^{\infty} \gamma_j (\sigma - 1)^j \;, &(2.13)
\cr}
$$
where $\zeta(s)$ is Riemann's zeta function, each of $\gamma$,
$\gamma_1$, $\gamma_2$, ...  is a certain (real valued) absolute
constant (the first of these, $\gamma$, being Euler's constant)
and $\theta(K , \sigma)$ is a number lying in the interval~$(0,1)$.
By~(2.12), we have $\alpha = -\zeta(1/2)$ in~(2.10), and we can calculate that
$\alpha - 1 = - (\zeta(1/2) + 1) = 0.4603545\ldots\ $. Given that
$\zeta(2) = \pi^2 /6$, we find (similarly) that
$\zeta_2 = (\pi^2 /6) - N^{-1} + O(N^{-2})$
in (2.3) to~(2.7). We also note that $\zeta_1 = \log N  + \gamma + O(1/N)$
(as follows, for example, by letting $\sigma\rightarrow 1$ in~(2.13)).

\par
We remark that, by combining methods similar to those used to obtain
(2.10) with certain applications of the Euler-Maclaurin
summation formula, we have been able to determine that the
variable $\phi\in [0,1)$ in (2.3) and~(2.7) satisfies
$$
\phi = \beta + O\left( {1 + \log N \over N^{1/7}}\right)\;,
\eqno(2.14)
$$
where
$\beta = 1 - {\pi^2 \over 24} - {1\over 2}(\log(2\pi) - 1)^2
+ {1\over 2}(1-\gamma)^2 =0.32712\ldots\ $.
We omit our proof of~(2.14), which shows no features that are truly
novel (and would require more than just a few pages). By~(2.14),
we can sharpen (2.8) somewhat, for large values of~$N$.
\par
Finally we consider the estimation problem~(d), stated earlier.
Noting firstly that ${\bf w}={\bf u} - (\zeta_1 /\zeta_2){\bf f}$, we
are able to deduce that
$$
\| {\bf w}\|^2 = N - {\zeta_1^2 \over \zeta_2}
=N + O\left( (1 + \log N)^2\right)
\eqno(2.15)
$$
and that
$$
{\bf w}^T A {\bf w}
={\bf u}^T A {\bf u} - 2\left( \zeta_1 /\zeta_2\right) {\bf u}^T A {\bf f}
+ \left( \zeta_1 /\zeta_2\right)^2 {\bf f}^T A {\bf f}\;.
\eqno(2.16)
$$
We have, moreover,
$$
{\bf u}^T A {\bf u}
= \sum_{m\leq N} \sum_{n\leq N} \left[ {N^2 \over mn}\right]
= \sum_m \sum_n \left[ {N^2 \over mn}\right]
  -2 \sum_{m > N} \sum_n \left[ {N^2 \over mn}\right]
= D_1 - 2 D_2\quad\ \hbox{(say).}
\eqno(2.17)
$$
Here
$$
D_1=\sum_{\ell\leq N^2} \tau_3(\ell)
=\left( {1\over 2}\log^2\left( N^2\right) +(3\gamma -1)\log\left( N^2\right)
+ c_1\right) N^2 + O\left( N^{\varepsilon + 43/48}\right)\;,
\eqno(2.18)
$$
where $c_1 = 3\gamma^2 - 3\gamma + 3\gamma_1 + 1$;
see pages 352-4 of~[{\bf 4}]
for the second equality in~(2.18).
\par
Regarding the sum $D_2$ in~(2.17), we have:
$$
\eqalign{
D_2 =\ \sum_{m>N}\ \sum
\sum_{\!\!\!\!\!\!\!\!\!\!\!\!\!\!
{\scriptstyle n\quad\  k\atop\scriptstyle (nk)m\leq N^2}} 1
 &=\sum_{\ell < N} \biggl( \sum_{n \mid \ell} 1\biggr)
 \sum_{N < m \leq N^2 /\ell} 1 \cr
 &= N^2 \sum_{\ell < N} {\tau_2 (\ell)\over\ell}
 - N \sum_{\ell < N} \tau_2 (\ell)
 + O\left( \sum_{\ell < N} \tau_2 (\ell)\right) \;.
\cr}
$$
By partial summation and Huxley's estimate
on page~593 of~[{\bf 3}]
for the remainder term in Dirichlet's divisor problem
(namely $\Delta(x)=\sum_{\ell\leq x} \tau_2(\ell) - (\log x + 2\gamma -1)x$),
we deduce from the above that
$$
D_2 = \left( {1\over 2}\log^2 N + (2\gamma -1)\log N + c_2\right) N^2
+ O\left( N^{547/416} (\log N)^{3.26}\right) \;,
$$
where
$$
c_2 = \int_1^{\infty} {\Delta(x) {\rm d}x\over x^2}
=\lim_{\sigma\rightarrow 2+} \left( {\zeta^2 (\sigma -1)\over \sigma -1}
- {1\over (\sigma -2)^2} - {2\gamma -1\over \sigma - 2}\right)
= \gamma^2 - 2\gamma + 2\gamma_1 +1
$$
(with $\gamma$ and $\gamma_1$ as in~(2.13)).
Using this,(2.17), and~(2.8), we have
$$
{\bf u}^T A {\bf u}
= \left( \log^2 N + 2\gamma \log N + c_3\right) N^2
+O\left( N^{547/416} (\log N)^{3.26}\right) \;,
\eqno(2.19)
$$
where $c_3=c_1 - 2c_2 = \gamma^2 + \gamma - \gamma_1 -1$.
Trivial estimates show that one has
${\bf u}^T A {\bf f} = \zeta_1 \zeta_2 N^2 +O((1+\log N)N)$:
using this, (2.4), (2.19), (2.16), and estimates already obtained for
$\zeta_1$ and~$\zeta_2$, we find that
$$
\eqalignno{
{\bf w}^T A {\bf w}
 &= \left( \log^2 N + 2\gamma \log N + c_3 - \zeta_1^2\right) N^2
+O\left( N^{547/416} (\log N)^{3.26}\right) \cr
 &= c_4 N^2 + O\left( N^{547/416} (\log N)^{3.26}\right)\;, &(2.20)
\cr}
$$
where
$c_4 = c_3 - \gamma^2 =
\gamma - \gamma_1 -1 = 0.57721566\ldots\ - 0.07281584\ldots\ -1
= -0.495600\ldots\ $
(see~[{\bf 10}].

\par
Since (2.15) implies $N > \| {\bf w}\|^2 \geq N/10$, then for $N\geq 2$,
using (2.15), (2.20), and Rayleigh's principle shows that:
$$
\lambda_1 \leq {{\bf w}^T A {\bf w} \over \| {\bf w}\|^2}
< c_4 N + O\left( N^{131/416} (\log N)^{3.26}\right)
\quad\ \hbox{($N\geq 2$).}
\eqno(2.21)
$$
The coefficient of $N$ in this upper bound may well be close to optimal:
when $N=10321$, for example,
computations done with the `GNU Octave' software package returned
$-0.493678\ldots\ $
as an estimate of the value of $\lambda_1 /N$ in this case.
By reasoning similar to that which gives~(2.9), we may deduce
from (2.7), (2.14) and~(2.21) that, as $N\rightarrow\infty$, we have
$|\lambda_2|/N < (1+o(1)) (\beta - c_4^2)^{1/2} \sim 0.2855539\ldots\ $
and
$({\bf e}_1 \cdot \hat{\bf w})^2
\geq (0.5 + o(1))(1 + (2 c_4^2 \beta^{-1} - 1)^{1/2})
\sim 0.8540699\ldots\ $. Therefore, for $N$ sufficiently large,
the lines $\{ t{\bf w} \,:\, t\in{\Bbb R}\}$ and
$\{ t{\bf e}_1 \,:\, t\in{\Bbb R}\}$
will meet at an angle of less than $\pi /8$ radians.
\par
We end this section with some speculations driven by certain numerical
evidence, gathered with the help of `GNU Octave'. We omit the detailed
evidence, but instead just summarise what it suggests.
Let $k$ be any fixed non-zero integer, and let $N$ now be free to vary
in the range $N > |k|$. Our numerical evidence suggests that
$\lambda_{\left\{ -k/N\right\} N}\sim \Lambda_k N$ as $N\rightarrow\infty$,
where $\Lambda_k$ is a real number that depends only on~$k$, and
where each of the two associated sequences,
$\Lambda_1 , \Lambda_2 , \Lambda_3 , \ldots\ $ and
$-\Lambda_{-1} , -\Lambda_{-2} , -\Lambda_{-3} , \ldots\ $, decreases
monotonically, and converges to~$0$. Further numerical evidence
suggests that if  $\theta\in (0,1)$ is fixed, and if
$e_{j,\ell}$ denotes the $\ell$-th component of the normalised
eigenvector~${\bf e}_j$, so that
${\bf e}_j = ( e_{j,1} , e_{j,2} , \ldots , e_{j,N} )^{\rm T}$ for
$j = 1,\ldots ,N$,
then as $N\rightarrow \infty$ we appear to see that
$$
e_{\left\{ -k/N\right\} N , \ell}
\sim(-1)^{b(N,k)}N^{-1/2} E_k (\ell /N)
\quad\hbox{  for $\ell = [\theta N]+1, [\theta N]+2, \ldots , N$,}
$$
with $E_k$ here being a certain real function independent of $\ell$ and~$N$
that is continuous on~$(0,1]$,
and with an integer exponent $b(N,k)$ independent of~$\ell$.
The occurrence of the functions $E_{\pm 1} , E_{\pm 2} , E_{\pm 3} , \ldots\ $
in this might be explained if they were eigenfunctions of a suitable
linear operator ${\cal A} : L^2 [0,1]\rightarrow L^2 [0,1]$.

\beginsection 3.
Various Decompositions of ${\bf m}^T A {\bf m}$ in the principal case

\medskip

It is our hope (as yet unrealised) that a study of the quadratic form
${\bf v}^T A {\bf v}$ (particularly when ${\bf v}$ is the vector
${\bf m}=(\mu(1),\ldots ,\mu(N))^T$), in the principal case of~(1.5),
might lead to new results about the Mertens function $M(x)$. In this
section we briefly describe (and compare) several different approaches
to such an investigation, each involving a different decomposition of
the quadratic form. We find it convenient to modify the earlier
notation $M(g,x)$ in (1.3): we use $M(s, x)$, where $s$ is a complex
number, (rather than a function), to mean $M(g, x)$ for the power
function $g(n) = n^{-s}$.
\par
We consider firstly (2.1) with ${\bf v} = {\bf m}$. We assume throughout
that $N$ is large. As the eigenvalue $\lambda_N$ is exceptionally large
among all the eigenvalues of~$A$,
we handle the term $\lambda_N ({\bf e}_N \cdot {\bf m})^2$ with some care.
As substitution of $-{\bf e}_N$ for~${\bf e}_N$ does not alter this
term, we can take the ambiguous sign in~(2.9) to be~$+$. We note that
$$
({\bf e}_N \cdot {\bf m})^2
= \left( ({\bf e}_N - \hat{\bf f}) \cdot {\bf m}\right)^2
  + 2 \left( ({\bf e}_N
  - \hat{\bf f}) \cdot {\bf m}\right) (\hat{\bf f}\cdot {\bf m})^2
  + (\hat{\bf f} \cdot {\bf m})^2 \;.
\eqno(3.1)
$$
Here
$$
\hat{\bf f} \cdot {\bf m}
=  \| {\bf f}\|^{-1} {\bf f} \cdot {\bf m}
=  {1\over\sqrt{\zeta_2}} \sum_{n\leq N} {\mu(n)\over n}
=  {M(1,N)\over\sqrt{\zeta_2}} \ll \log N
$$
and, by the Cauchy-Schwarz inequality and~(2.9),
$$
| ({\bf e}_N - \hat{\bf f}) \cdot {\bf m}|
  \leq \| {\bf e}_N - \hat{\bf f}\| \cdot \| {\bf m}\|
=  O\left( {\log N\over N}\right) \cdot \sqrt{\sum_{n\leq N} \mu^2(n)}
  \ll {\log N\over\sqrt{N}}\;.
$$
By these results, together with (3.1) and~(2.6), we have:
$$
\lambda_N ( {\bf e}_N \cdot {\bf m} )^2
= O\left( N\log^2 N\right)
  + O\left( N^{3/2} (\log N) |M(1,N)| \right) + N^2 (M(1,N))^2 \;.
\eqno(3.2)
$$
\par
Small eigenvalues make a relatively insignificant contribution here, for
(2.2) and (2.8) imply that if $1\leq K\leq N/2$, then
$$
\sum_{k=K}^{N-K} |\lambda_k| \left( {\bf e}_k \cdot {\bf m}\right)^2
< {N\over\sqrt{K}} \sum_{n=1}^N \left( {\bf e}_k \cdot {\bf m}\right)^2
= {N\over\sqrt{K}} \| {\bf m} \|^2 \leq {N^2 \over\sqrt{K}} \;.
$$
By this, and by (3.2) and (2.1) (for ${\bf v} = {\bf m}$), we find that
$$
\eqalign{
{{\bf m}^T A {\bf m} \over N^2}
 &= (M(1,N))^2\, +\ \left( \| {\bf m}\|^2 / N\right)  \!\!\!\!\!\!\!\!
 \sum_{\scriptstyle 1 \leq k < N\atop\scriptstyle \min\{ k , N-k\} < K}
  \!\!\!\!\!\!\!\!\left( \lambda_k / N\right)
  \left( {\bf e}_k \cdot \hat{\bf m}\right)^2 \cr
 &\quad\ + O\left( K^{-1/2} + N^{-1/2} (\log N) |M(1,N)| +
 N^{-1}\log^2 N \right) \;,
\cr}
\eqno(3.3)
$$
for $K=1,2,\ldots , N^2$. We remark that, if the second of the three terms
on the right-hand side of~(3.3) is
considered in isolation, then we observe trivially from~(2.8) that
the absolute value of this term is $O(\sqrt{K}\,)$. Taking account of
the context here (the relation (3.3) and the principal case of~(1.5)
and~(1.3)), and noting also that $|M(1,N)|\leq \| {\bf m}\|^2 / N\,$
(a consequence of (1.11), the trivial bound $|[y] - y| < 1$, and the
fact that $[N/1]-(N/1) = 0$, it is clear that this term is a bounded
function of the pair $(N,K)\in{\Bbb N}^2$. This gives some idea of the
gap that must be bridged if (3.3) is to help in the study of~$M(x)$.
\par
To reach (3.3) we have used the work of Section~2, on $\lambda_N$ and
${\bf e}_N$.
Our next decomposition of ${\bf m}^T A {\bf m}$ avoids such results, but
nevertheless has much in common with~(3.3).
\par
First we use $[x] = x - {1\over 2} - \psi(x)$, where
$\psi(x) = \{ x\} - {1\over 2}$. We have
$$
A = N^{2\,} {\bf f}^{\,} \!{\bf f}^T -
\textstyle{1\over 2} {\bf u} {\bf u}^T + Z\;,
\eqno(3.4)
$$
where $Z$ is the $N\times N$ matrix  of elements
$z_{mn} = - \psi(N^2 / (mn))$, whilst ${\bf f}$ and ${\bf u}$ are as in
Section~2. We have trivially ${\rm Tr}(Z^2) < N^2 /4$; with the help of
(2.14), (2.19), and an estimate for~$\zeta_1$, we obtain the sharper
result that ${\rm Tr}(Z^2) \sim c_5 N^2$ as $N\rightarrow\infty$,
where $c_5 = \beta + {1\over 4} + c_3 - \gamma^2 = 0.0815206\ldots\ $.
Reasoning as in the derivation of~(3.3), we see from~(3.4) that, for
$K=1,2,\ldots , N^2\,$ (say), one has
$$
\eqalignno{
{{\bf m}^T A {\bf m} \over N^2}
 &= \left( {\bf m} \cdot {\bf f}\right)^2
 - {\left( {\bf m} \cdot {\bf u}\right)^2 \over 2 N^2}
 + {{\bf m}^T Z {\bf m} \over N^2} &(3.5) \cr
 &= (M(1,N))^2\, -{(M(N))^2\over 2 N^2} \cr
 &\quad\ +\ \left( \| {\bf m}\|^2 / N\right)  \!\!\!\!\!\!\!\!
 \sum_{\scriptstyle 1 \leq k \leq N\atop\scriptstyle \min\{ k , N+1-k\} < K}
 \!\!\!\!\!\!\!\!\left( \widetilde\lambda_k / N\right)
 \left( \widetilde{\bf e}_k \cdot \hat{\bf m}\right)^2
 + O\left( K^{-1/2} \right)\;, &(3.6)
\cr}
$$
where
$\widetilde\lambda_1 \leq \widetilde\lambda_2 \leq\ \cdots \
\leq \widetilde\lambda_N$ are the eigenvalues of $Z$,
while $\widetilde{\bf e}_1 , \ldots , \widetilde{\bf e}_N$
form the corresponding orthonormal basis of eigenvectors.
We note the presence of the term $-{1\over 2}N^{-2}(M(N))^2$ in~(3.6),
which is not apparent in~(3.3): in view of our results on Problem~(d)
of Section~2, one may regard this term as being an approximation to the
term
$(\| {\bf m}\|^2 / N) ( \lambda_1 / N) ( {\bf e}_1 \cdot \hat{\bf m})^2
 = N^{-2} \lambda_1 ( {\bf e}_1 \cdot {\bf m})^2$, which is present
in~(3.3) for $K > 1$.
\par
We remark that (3.5) permits an alternative, non-spectral, decomposition of
${\bf m}^T A {\bf m}$, through substituting the usual truncated Fourier
expansion of the function $\psi $ into each element of the matrix~$Z$
in~(3.5):
$$
-\psi(x) = \sum_{0<h\leq H} {\sin(2\pi h x) \over \pi h}
  + O\left( {\eta\over \eta
  + \min\{ |x-\ell | \,:\, \ell\in{\Bbb Z}\}}\right)
  \quad\ \hbox{($H=1/\eta\geq 1$).}
$$
This leads (via estimates
from~[{\bf 11}])
to the decompositions
$$
{\bf m}^T Z {\bf m} = \sum_{h=1}^H {{\bf m}^T Z(h) {\bf m} \over \pi h}
 + O\left( {N^2 (\log N)^2 \log H \over H}\right)
 \quad\ \hbox{(for $H=1,2,\ldots , N$ (say)),}
$$
where $Z(h)$ is the $N\times N$ matrix with elements
$z_{mn}(h)=\sin(2\pi h N^2 /(mn))$.
We have yet to explore making proper use of this truncation idea.
\par
One further approach to the decomposition of ${\bf m}^T A {\bf m}$ uses
Perron's formula, Theorem~5.1 of~[{\bf 9}],
equation~(A.8) of~[{\bf 4}].
We apply Perron's formula as in
Lemma~3.12 of~[{\bf 12}],
adapting the proof to sharpen certain error terms (parts of the
improvement are results
of Shiu~[{\bf 11}]).
We find that if, whenever ${\rm Re}(s) > 1$, one has
$$
F(s) = \sum_{\ell =1}^{\infty} {a_{\ell}\over \ell^s}
=  \left( \sum_{m\leq y} {\alpha_m \over m^s}\right)
  \left( \sum_{n\leq z} {\beta_n \over n^s}\right) \zeta(s)
=  A(s) B(s) \zeta(s) \quad\ \hbox{(say),}
\eqno(3.7)
$$
where $y,z \geq 1$ and $\alpha_m,\beta_n$ denote complex constants of
modulus less than or equal to $1$, then, for any fixed $\varepsilon > 0$,
when $x = yz$, in the ranges $1 < c \leq 2$ and
$3\leq T \leq x^{1-\varepsilon}$, we have
$$
{1\over 2\pi i}\int\limits_{c-iT}^{c+iT} F(s) x^s \, {{\rm d}s\over s}
= \sum_{\ell\leq x} a_{\ell}
+ O\left( {x^c \log^2 x \over (c-1) T}\right)
+ O_{\varepsilon}\left( {x (\log x)^2 (\log T) \over T}\right) \;.
\eqno(3.8)
$$
To link this to our matrix~$A$, we observe that (3.7) implies
$$
\sum_{\ell\leq x} a_{\ell}
=  \sum_{\ell\leq x}\ \sum\ \sum\  \sum_{
  \!\!\!\!\!\!\!\!\!\!\!\!\!\!\!\!\!\!\!\!\!\!\!\!\!\!\!\!\!\!\!\!\!
  {\scriptstyle m\leq y\ n\leq z\quad k\atop\scriptstyle
  \quad  mnk = \ell}} \alpha_m \beta_n
= \ \sum\ \sum\  \sum_{
  \!\!\!\!\!\!\!\!\!\!\!\!\!\!\!\!\!\!\!\!\!\!\!\!\!\!\!\!\!\!\!\!\!
  {\scriptstyle m\leq y\ n\leq z\quad k\atop\scriptstyle
  \quad\  mnk \leq x}} \alpha_m \beta_n
=  \sum_{m\leq y}\,\sum_{n\leq z} \left[ {x\over mn}\right]
  \alpha_m \beta_n\;.
$$
\par
Setting $c=1 + (\log x)^{-1}$ in (3.8), we shift the contour of
integration there until it aligns with the line ${\rm Re}(s) = {1\over 2}$:
in so doing, we pick up a contribution from the residue of $\zeta(s)$
at its pole, $s = 1$, and also some remainder terms, which are integrals
along the line segments joining ${1\over 2}+iT$ to $c+iT$, and
${1\over 2}-iT$ to $c-iT$.
By Theorem~7.2~(A) of Titchmarsh~[{\bf 12}],
we deduce that these remainder term integrals are of size
$O(x (\log x)^2 \sqrt{\log T} /T)$ for almost all values of~$T$ (in a
certain sense) lying in any given `dyadic interval'
$[T_0 , 2T_0]\subseteq [3 , 2 x^{1-\varepsilon}]$.
Hence we arrive at the conclusion that, for any given $\varepsilon > 0$
when $x = yz$ and $3\leq T_0\leq x^{1-\varepsilon}$, we have
$$
\sum_{m\leq y}\,\sum_{n\leq z} \left[ {x\over mn}\right] \alpha_m \beta_n
=  A(1)B(1)x
+ {1\over 2\pi i}\int\limits_{{1\over 2} - iT}^{{1\over 2} + iT}
  A(s)B(s)\zeta(s) x^s \, {{\rm d}s\over s}
 + O_{\varepsilon}\left( {x\log^3 x\over T}\right) \;,
$$
for some $T\in [T_0 , 2T_0]$. We specialise this to the case
$\varepsilon = 1/2$, $y = z = N$, where $N$ is a positive integer, so
that $x = N^2$, and $\alpha_n = \beta_n = \mu(n)\,$. We find that
when $3\leq T_0 \leq N$, there exists some $T\in [T_0 , 2 T_0]$ such that
$$
\eqalign{
{{\bf m}^T A {\bf m} \over N^2}
 &= (M(1,N))^2
  +{\| {\bf m}\|^2 \over N} \int_{-T}^T {\zeta_1 N^{2it}
  \zeta\left( {1\over 2} + it\right) \over \left( \pi + 2\pi it\right)}
  \left( {M\left( {1\over 2} + it , N\right) \over \sqrt{\zeta_1}
  \| {\bf m}\|}\right)^2 {\rm d}t \cr
 &\quad\ +O\left( T_0^{-1} \log^3 N\right) \;.
\cr}
\eqno(3.9)
$$
If we put ${\bf E}(s)=(1^{-s} , 2^{-s} , \ldots , N^{-s})^T\,$ for a
fixed complex number~$s$,
then the factor $M({1\over 2} + it , N) / (\sqrt{\zeta_1} \| {\bf m}\|)$
here may be written as $\hat{\bf E}({1\over 2} + it) \cdot \hat{\bf m}$:
the decomposition in~(3.9) may therefore be considered similar in form
to that in~(3.3), although (3.9) involves an integration over the
range~$[-T,T]$, instead of the summation over a subset of the
(discrete) spectrum of $A$ that we had in~(3.3).

\beginsection References

\item{[{\bf 1}]} J.-P. Cardinal, `Symmetric matrices related to the Mertens
function',  {\it Linear Algebra Appl.} {\bf 432} (2010), 161-172.
\item{[{\bf 2}]} D.R. Heath-Brown, `Prime Numbers in Short Intervals and a
Generalised Vaughan Identity',
{\it Can. J. Math.}, {\bf 34}, no. 6 (1982), 1365-1377.
\item{[{\bf 3}]} M.N. Huxley, `Exponential Sums and Lattice Points III',
{\it Proc. London Math. Soc.} (3), {\bf 87} (2003), 591-609.
\item{[{\bf 4}]} A. Ivi\'{c}, {\it The Riemann Zeta-Function}, Dover
Publications, Mineola, New York (2003).
\item{[{\bf 5}]} H. Iwaniec and E. Kowalski, {\it Analytic Number
Theory}, A.M.S. Colloquium Publications 53,  American Mathematical
Society, Providence RI, 2004.
\item{[{\bf 6}]} Yu. V. Linnik, `All large numbers are sums of a prime
and two squares I', {\it Mat. Sbornik Nov. Ser.} {\bf 52 (94)} (1960),
561-700.
\item{[{\bf 7}]} Yu. V. Linnik, `All large numbers are sums of a prime
and two squares II', {\it Mat. Sbornik Nov. Ser.} {\bf 53 (95)} (1961),
3-38.
\item{[{\bf 8}]} E. Meissel, `Observationes quaedam in theoria numerorum',
{\it J. Reine Angew. Math.}, {\bf 48} (1854), 301-316.
\item{[{\bf 9}]} H.L. Montgomery and R.C. Vaughan, {\it Multiplicative
Number Theory I. Classical Theory},  Cambridge Studies in Advanced
Mathematics 97, Cambridge University Press (2007).
\item{[{\bf 10}]} {\it On-Line Encyclopedia of Integer Sequences}
(``Sloane's''), http://oeis.org
\item{[{\bf 11}]} P. Shiu, `A Brun-Titchmarsh theorem for
multiplicative functions',
{\it J. Reine Angew. Math.}, {\bf 313} (1980), 161-170.
\item{[{\bf 12}]} E.C. Titchmarsh (revised by D.R. Heath-Brown), {\it The
Theory of the Riemann Zeta-function}, Oxford Univ. Press, 1986
\item{[{\bf 13}]} R.C. Vaughan, `An Elementary Method in Prime Number Theory',
{\it Recent Progress in Analytic Number Theory, vol.} 1 (Durham, 1979),
Academic Press, London - New York, 1981, pp. 341-348.

\bye